\documentclass[reqno,12pt,a4paper]{amsart}

\usepackage{mathrsfs}
\usepackage{amssymb}
\usepackage{amsfonts}
\usepackage{latexsym}
\usepackage{amsthm}

\usepackage{graphicx}
\def\lb{\label}

\newcommand{\er}[1]{\textrm{(\ref{#1})}}

\begin{document}


\renewcommand{\theequation}{\arabic{section}.\arabic{equation}}
\theoremstyle{plain}
\newtheorem{theorem}{\bf Theorem}[section]
\newtheorem{lemma}[theorem]{\bf Lemma}
\newtheorem{corollary}[theorem]{\bf Corollary}
\newtheorem{proposition}[theorem]{\bf Proposition}
\newtheorem{definition}[theorem]{\bf Definition}
\newtheorem{remark}[theorem]{\it Remark}

\def\a{\alpha}  \def\cA{{\mathcal A}}     \def\bA{{\bf A}}  \def\mA{{\mathscr A}}
\def\b{\beta}   \def\cB{{\mathcal B}}     \def\bB{{\bf B}}  \def\mB{{\mathscr B}}
\def\g{\gamma}  \def\cC{{\mathcal C}}     \def\bC{{\bf C}}  \def\mC{{\mathscr C}}
\def\G{\Gamma}  \def\cD{{\mathcal D}}     \def\bD{{\bf D}}  \def\mD{{\mathscr D}}
\def\d{\delta}  \def\cE{{\mathcal E}}     \def\bE{{\bf E}}  \def\mE{{\mathscr E}}
\def\D{\Delta}  \def\cF{{\mathcal F}}     \def\bF{{\bf F}}  \def\mF{{\mathscr F}}
\def\c{\chi}    \def\cG{{\mathcal G}}     \def\bG{{\bf G}}  \def\mG{{\mathscr G}}
\def\z{\zeta}   \def\cH{{\mathcal H}}     \def\bH{{\bf H}}  \def\mH{{\mathscr H}}
\def\e{\eta}    \def\cI{{\mathcal I}}     \def\bI{{\bf I}}  \def\mI{{\mathscr I}}
\def\p{\psi}    \def\cJ{{\mathcal J}}     \def\bJ{{\bf J}}  \def\mJ{{\mathscr J}}
\def\vT{\Theta} \def\cK{{\mathcal K}}     \def\bK{{\bf K}}  \def\mK{{\mathscr K}}
\def\k{\kappa}  \def\cL{{\mathcal L}}     \def\bL{{\bf L}}  \def\mL{{\mathscr L}}
\def\l{\lambda} \def\cM{{\mathcal M}}     \def\bM{{\bf M}}  \def\mM{{\mathscr M}}
\def\L{\Lambda} \def\cN{{\mathcal N}}     \def\bN{{\bf N}}  \def\mN{{\mathscr N}}
\def\m{\mu}     \def\cO{{\mathcal O}}     \def\bO{{\bf O}}  \def\mO{{\mathscr O}}
\def\n{\nu}     \def\cP{{\mathcal P}}     \def\bP{{\bf P}}  \def\mP{{\mathscr P}}
\def\r{\varrho} \def\cQ{{\mathcal Q}}     \def\bQ{{\bf Q}}  \def\mQ{{\mathscr Q}}
\def\s{\sigma}  \def\cR{{\mathcal R}}     \def\bR{{\bf R}}  \def\mR{{\mathscr R}}
\def\S{\Sigma}  \def\cS{{\mathcal S}}     \def\bS{{\bf S}}  \def\mS{{\mathscr S}}
\def\t{\tau}    \def\cT{{\mathcal T}}     \def\bT{{\bf T}}  \def\mT{{\mathscr T}}
\def\f{\phi}    \def\cU{{\mathcal U}}     \def\bU{{\bf U}}  \def\mU{{\mathscr U}}
\def\F{\Phi}    \def\cV{{\mathcal V}}     \def\bV{{\bf V}}  \def\mV{{\mathscr V}}
\def\P{\Psi}    \def\cW{{\mathcal W}}     \def\bW{{\bf W}}  \def\mW{{\mathscr W}}
\def\o{\omega}  \def\cX{{\mathcal X}}     \def\bX{{\bf X}}  \def\mX{{\mathscr X}}
\def\x{\xi}     \def\cY{{\mathcal Y}}     \def\bY{{\bf Y}}  \def\mY{{\mathscr Y}}
\def\X{\Xi}     \def\cZ{{\mathcal Z}}     \def\bZ{{\bf Z}}  \def\mZ{{\mathscr Z}}
\def\O{\Omega}

\newcommand{\mc}{\mathscr {c}}

\newcommand{\gA}{\mathfrak{A}}          \newcommand{\ga}{\mathfrak{a}}
\newcommand{\gB}{\mathfrak{B}}          \newcommand{\gb}{\mathfrak{b}}
\newcommand{\gC}{\mathfrak{C}}          \newcommand{\gc}{\mathfrak{c}}
\newcommand{\gD}{\mathfrak{D}}          \newcommand{\gd}{\mathfrak{d}}
\newcommand{\gE}{\mathfrak{E}}
\newcommand{\gF}{\mathfrak{F}}           \newcommand{\gf}{\mathfrak{f}}
\newcommand{\gG}{\mathfrak{G}}           
\newcommand{\gH}{\mathfrak{H}}           \newcommand{\gh}{\mathfrak{h}}
\newcommand{\gI}{\mathfrak{I}}           \newcommand{\gi}{\mathfrak{i}}
\newcommand{\gJ}{\mathfrak{J}}           \newcommand{\gj}{\mathfrak{j}}
\newcommand{\gK}{\mathfrak{K}}            \newcommand{\gk}{\mathfrak{k}}
\newcommand{\gL}{\mathfrak{L}}            \newcommand{\gl}{\mathfrak{l}}
\newcommand{\gM}{\mathfrak{M}}            \newcommand{\gm}{\mathfrak{m}}
\newcommand{\gN}{\mathfrak{N}}            \newcommand{\gn}{\mathfrak{n}}
\newcommand{\gO}{\mathfrak{O}}
\newcommand{\gP}{\mathfrak{P}}             \newcommand{\gp}{\mathfrak{p}}
\newcommand{\gQ}{\mathfrak{Q}}             \newcommand{\gq}{\mathfrak{q}}
\newcommand{\gR}{\mathfrak{R}}             \newcommand{\gr}{\mathfrak{r}}
\newcommand{\gS}{\mathfrak{S}}              \newcommand{\gs}{\mathfrak{s}}
\newcommand{\gT}{\mathfrak{T}}             \newcommand{\gt}{\mathfrak{t}}
\newcommand{\gU}{\mathfrak{U}}             \newcommand{\gu}{\mathfrak{u}}
\newcommand{\gV}{\mathfrak{V}}             \newcommand{\gv}{\mathfrak{v}}
\newcommand{\gW}{\mathfrak{W}}             \newcommand{\gw}{\mathfrak{w}}
\newcommand{\gX}{\mathfrak{X}}               \newcommand{\gx}{\mathfrak{x}}
\newcommand{\gY}{\mathfrak{Y}}              \newcommand{\gy}{\mathfrak{y}}
\newcommand{\gZ}{\mathfrak{Z}}             \newcommand{\gz}{\mathfrak{z}}

\def\ve{\varepsilon}   \def\vt{\vartheta}    \def\vp{\varphi}    \def\vk{\varkappa}

\def\A{{\mathbb A}} \def\B{{\mathbb B}} \def\C{{\mathbb C}}
\def\dD{{\mathbb D}} \def\E{{\mathbb E}} \def\dF{{\mathbb F}} \def\dG{{\mathbb G}}
\def\H{{\mathbb H}}\def\I{{\mathbb I}} \def\J{{\mathbb J}} \def\K{{\mathbb K}}
\def\dL{{\mathbb L}}\def\M{{\mathbb M}} \def\N{{\mathbb N}} \def\dO{{\mathbb O}}
\def\dP{{\mathbb P}} \def\R{{\mathbb R}} \def\dQ{{\mathbb Q}}
\def\S{{\mathbb S}} \def\T{{\mathbb T}} \def\U{{\mathbb U}}
\def\V{{\mathbb V}}\def\W{{\mathbb W}} \def\X{{\mathbb X}} \def\Y{{\mathbb Y}} \def\Z{{\mathbb Z}}

\newcommand{\1}{\mathbbm 1}
\newcommand{\dd}    {\, \mathrm d}



\def\la{\leftarrow}              \def\ra{\rightarrow}            \def\Ra{\Rightarrow}
\def\ua{\uparrow}                \def\da{\downarrow}
\def\lra{\leftrightarrow}        \def\Lra{\Leftrightarrow}


\def\lt{\biggl}                  \def\rt{\biggr}
\def\ol{\overline}               \def\wt{\widetilde}
\def\no{\noindent}


\let\ge\geqslant                 \let\le\leqslant
\def\lan{\langle}                \def\ran{\rangle}
\def\/{\over}                    \def\iy{\infty}
\def\sm{\setminus}               \def\es{\emptyset}
\def\ss{\subset}                 \def\ts{\times}
\def\pa{\partial}                \def\os{\oplus}
\def\om{\ominus}                 \def\ev{\equiv}
\def\iint{\int\!\!\!\int}        \def\iintt{\mathop{\int\!\!\int\!\!\dots\!\!\int}\limits}
\def\el2{\ell^{\,2}}             \def\1{1\!\!1}
\def\sh{\sharp}
\def\wh{\widehat}

\def\all{\mathop{\mathrm{all}}\nolimits}
\def\where{\mathop{\mathrm{where}}\nolimits}
\def\as{\mathop{\mathrm{as}}\nolimits}
\def\Area{\mathop{\mathrm{Area}}\nolimits}
\def\arg{\mathop{\mathrm{arg}}\nolimits}
\def\const{\mathop{\mathrm{const}}\nolimits}
\def\det{\mathop{\mathrm{det}}\nolimits}
\def\diag{\mathop{\mathrm{diag}}\nolimits}
\def\diam{\mathop{\mathrm{diam}}\nolimits}
\def\dim{\mathop{\mathrm{dim}}\nolimits}
\def\dist{\mathop{\mathrm{dist}}\nolimits}
\def\Im{\mathop{\mathrm{Im}}\nolimits}
\def\Iso{\mathop{\mathrm{Iso}}\nolimits}
\def\Ker{\mathop{\mathrm{Ker}}\nolimits}
\def\Lip{\mathop{\mathrm{Lip}}\nolimits}
\def\rank{\mathop{\mathrm{rank}}\limits}
\def\Ran{\mathop{\mathrm{Ran}}\nolimits}
\def\Re{\mathop{\mathrm{Re}}\nolimits}
\def\Res{\mathop{\mathrm{Res}}\nolimits}
\def\res{\mathop{\mathrm{res}}\limits}
\def\sign{\mathop{\mathrm{sign}}\nolimits}
\def\span{\mathop{\mathrm{span}}\nolimits}
\def\supp{\mathop{\mathrm{supp}}\nolimits}
\def\Tr{\mathop{\mathrm{Tr}}\nolimits}
\def\BBox{\hspace{1mm}\vrule height6pt width5.5pt depth0pt \hspace{6pt}}


\newcommand\nh[2]{\widehat{#1}\vphantom{#1}^{(#2)}}
\def\dia{\diamond}

\def\Oplus{\bigoplus\nolimits}




\def\qqq{\qquad}
\def\qq{\quad}
\let\ge\geqslant
\let\le\leqslant
\let\geq\geqslant
\let\leq\leqslant

\newcommand{\ca}{\begin{cases}}
\newcommand{\ac}{\end{cases}}
\newcommand{\ma}{\begin{pmatrix}}
\newcommand{\am}{\end{pmatrix}}
\renewcommand{\[}{\begin{equation}}
\renewcommand{\]}{\end{equation}}
\def\bu{\bullet}

\title[{Inverse  resonance scattering for on rotationally symmetric  manifolds}]
{Inverse  resonance scattering on rotationally symmetric  manifolds}

\date{\today}
\author[Hiroshi Isozaki]{Hiroshi Isozaki}
\address{Graduate School of Pure and Applied Sciences,
University of Tsukuba, Tsukuba, 305-8571, Japan,\ \
isozakih@math.tsukuba.ac.jp}
\author[Evgeny Korotyaev]{Evgeny Korotyaev }
\address{Department of Analysis,  Saint-Petersburg State University,   Universitetskaya nab. 7/9, St.
Petersburg, 199034, Russia, \ korotyaev@gmail.com, \
e.korotyaev@spbu.ru}


\begin{abstract}
\no  We discuss inverse resonance scattering for the Laplacian on a
rotationally symmetric manifold $M =  (0,\infty) \ts Y$ whose
rotation radius is constant outside  some compact interval. The
Laplacian on $M$ is unitarily equivalent to a direct sum of
one-dimensional Schr\"odinger operators with compactly supported
potentials on the half-line. We prove
\begin{itemize}
 \item Asymptotics of counting function of resonances at large radius.
\item
 Inverse problem: The rotation radius is uniquely determined by its
eigenvalues and resonances. Moreover, there exists an algorithm to
recover the rotation radius from its eigenvalues and resonances.
\end{itemize}
The
proof is based on some non-linear real analytic isomorphism between two Hilbert spaces.

\end{abstract}

\maketitle

\section {Introduction and main results}
\setcounter{equation}{0}

\subsection {Geometry.}
We consider  a rotationally symmetric
manifold $M = \R_+\ts Y$ equipped with a warped product
metric
\[
\lb{1} g = (dx)^2 + r^2(x)g_Y,\qq x\ge 0.
\]
Here $(Y, g_Y)$ is a compact m-dimensional Riemannian manifold (with
or without boundary), called tranversal manifold, and $r(x)>0$ is
the rotation radius satisfying
\begin{equation}
 r(x)=r_o=\const>0 \quad {\rm for} \quad  x\ge 1.
\label{r(x)geqr0fora>0}
\end{equation}
We introduce $q(x)$ by
\[
\lb{2}
\begin{aligned}
 r=r_o \,e^{{2\/m}Q},\qqq Q(x)=-\int_x^1 q(t)dt,\qq x\in \R_+=
 [0,\iy),
 \end{aligned}
\]
and assume that  the real function $q$ belongs to the class $\cP_1$, where $\cP_j$ is defined by
\[
\lb{dP} \cP_j =\left\{ g\, ; \,  g^{(j)}\in L^1_{real}(\R_+),\quad
\supp g\ss [0,1], \quad \sup \supp g=1  \right\}, \qq j\ge 0.
\]
As will be seen below, the geometry (i.e., the rotation
radius  $r$ and hence all derived quantities up to two integration
constants) is determined by $q$. The  Laplacian $\D_M$ on $M$ has the form
$$
- \D_M=-{1\/r^m}\pa_x \big(r^m\pa_x\big)-{\D_{Y}\/r^2},
$$
where $\D_{Y}$ is the Laplacian on $Y$. We assume the Dirichlet boundary condition, i.e. the domain of  $\D_{M} $ consists of the functions $f=f(x,y),(x,y)\in M=\R_+\ts Y$,
satisfying  the following boundary condition
\[
\lb{3u}
\begin{aligned}
 \ f|_{\pa M}=0,\qqq \pa
M=\{0\}
\end{aligned}
\]
The negative Laplacian $- \Delta_Y$ on $Y$ (with suitable boundary conditions when $Y$ has a boundary) has a discrete spectrum, $0\le E_1\le E_2\le  \cdots$,
 with corresonding orthonormal family
of eigenfunctions ${\P_\n, \n\ge 1}$, in $L^2(Y )$.

Let us reacll simple examples.

\noindent
(1) The case when $Y$ has no boundary. For example, if
$Y=\S^1$, then
$
E_1=0,
E_2=\pi^2, \cdots.$

\noindent
(2)  The case when $Y$ has a boundary. For example, if $Y$ is a compact interval in $\R^1$ e.g. $[0,1]$,
then
$
E_1=0,  E_2=\pi^2, \cdots,
$
for the Neumann boundary condition, and
$
E_1=\pi^2,
E_2=(2\pi)^2, \cdots,
$
for the Dirichlet boundary condition.

\medskip

Introduce the Hilbert spaces
$$
\cL_\n^2(M)=\rt\{h(x,y)=f(x)\P_\n(y)\, : \,  (x,y)\in
M=\R_+\ts Y, \qq \int_0^\iy|f(x)|^2dx<\iy
\rt\},\qq \n\ge 1.
$$
By the spectral decomposition of  $- \D_{Y}$, the Laplacian on $(M, g)$ acting in $L^2(M)=\os_{\n\ge
1} \cL_\n^2(M)$  is unitarily equivalent to a direct sum of
one-dimensional Schr\"odinger operators $H_\n$, namely,
\[
\lb{sum} - \D_{M}\backsimeq\os_{\n=1}^{\iy}\ (H_\n+u_{\n,0}).
\]
Here the  operator $H_\n$ acts on  $L^2(\R_+)$ under the
Dirichlet boundary condition and is given by
\[
\begin{aligned}
\lb{6k} H_\n f=-f''+p_\n f,\qq f(0)=0,
\end{aligned}
\]
and the potential $p_\n(x)$ is defined by
\[
\lb{dP1}
\begin{aligned}
\ca  & p_\n(x)=q'(x) +q^2(x)+u_\n(x)-u_{\n,0}, \qq  \\
&u_\n(x)= u_{\n,0}e^{-{4\/m}Q(x)},\qq Q(x)=-\int_x^1 q(t)dt,\
\\
&  u_{\n,0}={E_\n\/r^2(1)}.
\ac
\end{aligned}
\]
Note that
\begin{equation} \supp p_\n\ss[0,1], \qq p_{\nu} \in
L^1(0,1).
\nonumber
\end{equation}

 Below we fix an  index $\n$ arbitrarily, and omit it. We  consider the operator $H
f=-f''+pf$, where the potential $p=p_\n$ is given by \er{dP1}. It is
well-known (see \cite{F59}) that $H$ has purely absolutely
continuous spectrum $[0,\iy)$ and a finite number $N\ge 0$ of bound
states $E_1<...<E_{N}<0$ below the continuum.
The Jost solution $f_+(x,k )$ is a solution to the equation
\[
-f_+'' +p f_+=k ^2f_+,\ \ \ x\geq 0, \ \ \ k \in \C\sm \{0\}
\]
satisfying the condition $f_+(x,k )=e^{ixk },\ x\geq 1 $. The Jost
function $\psi(k)$ is defined by $\p(k)=f_+(0,k)$. Let $n_+(f)$ be
the number of zeros of $f$ analytic  in $\C_+$, each zero being
counted according to its multiplicity.
For an analytic function $f$ defined in a neighborhood of 0, we
define $n_o(f) = m$, if $f(z) =  z^mg(z)$ where $g(0) \neq 0$. The
Jost function $\psi(k)$ has the following standartd properties


$\bu $ The Jost function $\psi(k)$ is entire on $\C$ and has the
following asymptotics uniformly in $\arg k\in [0,\pi]$:
\[
\lb{Jp1} \p(k)=1+ {O(1/k)} \qqq \as \qq |k|\to \iy.
\]
$\bu $  The S-matrix is defined by
\[
S(k)={\p(-k)\/\p(k)}=e^{-i2\f_{sc}(k)},\ \ \ \ \f_{sc}(k)=\arg
\p(k), \ \ \ k \in \R .
\]
We call $\f_{sc}$  phase shift.
Note that $S(k)$ is meromorphic on $\C$. We define a resonance of $H$ to be a pole of $S(k)$ in $\C_-$.

 \noindent $\bu $  The zeros of $\p$ in $\C_+$  are given
by
\[
\lb{Jp2} k_{1}=i|E_1|^{1\/2},\dots k_{n_+}=i|E_{N}|^{1\/2}\in i\R_+,
\qqq N=n_+(\p).
\]
There may be a simple zero of $\psi$ at 0, hence $n_o(\psi) = 0$ or
1. The number of  zeros in $\C_-\cup \{0\}$ are infinite, and
arranged so that $0\leq |k_{{N}+1}|\leq |k_{N+2}|\leq \dots $
(counted with multiplicities). If $k_{{N}+1}=0$, then $0<|k_{N+2}|$.
The zeros of $\p$ in $\C_-$ coincide with the resonances.

\noindent $\bu $ Due to \er{Jp1}  we define the branch $\log \p(k),
k\in \ol\C_+\sm [0,k_1]$ by the condition $\log \p(iv)=o(1)$ as
$v\to \iy$. Since $\p(k)$ is real on $i\R_+$ and has $n_+$ zeros in
$\C_+$, the function $\f_{sc}$ is odd,  i.e. $ \f_{sc}(-k
)=-\f_{sc}(k )$ for  $k>0, $  continuous on $(0,\iy )$, $\f_{sc}(k
)=o(1)$ as $k\to \iy $ and
\[
\lb{Jp3} \f_{sc}(\pm 0)=\mp \pi \big (n_++{n_0\/2}\big),\ \ \ \
n_0=n_0(\p) = 0 \ {or} \ 1.
\]
Then, as a preliminary  result, we have

\begin{theorem}
\lb{T1} i) Let $q\in \cP_1$ and  $p$ be given by
\er{dP1}. Then $p\in \cP_0$.

\noindent
ii) If $r(x)\le r(1)$ for all $x\in [0,1]$, then the Laplacian
$-\D_M$ has no eigenvalues.

\no iii) There exists a radius  $r(x)>r_o$ for all $x\in (0,1)$
and  $\n_o$ large enough such
that each operator $\D_\n, \n>\n_o$ has eigenvalues and the
Laplacian $\D_M$ has infinite number of eigenvalues.

\end{theorem}

\vskip 10mm

\setlength{\unitlength}{2.0mm}
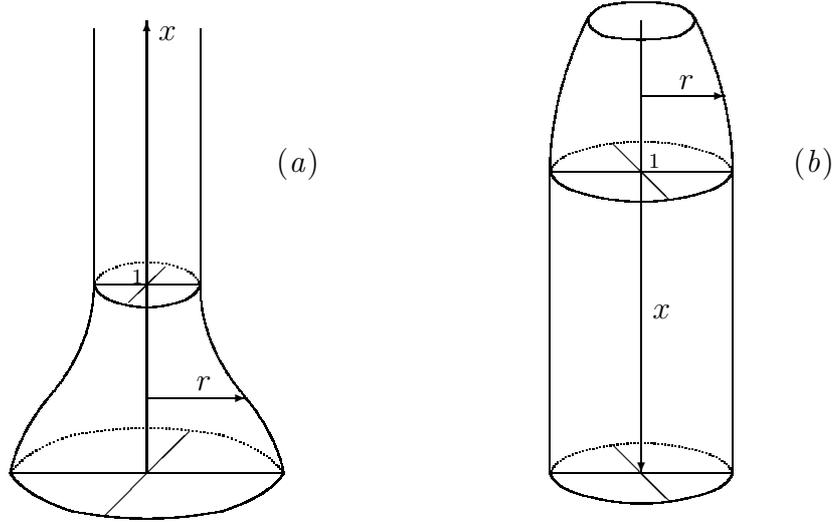
\begin{figure}[h]
\centering
\unitlength 1.0mm 
\begin{picture}(110,65)
\put(20,10){\vector(0,1){60.00}}

\bezier{150}(10,5)(20,3)(30,5)

\bezier{150}(10,5)(4,7)(2,10) \bezier{20}(10,15)(4,13)(2,10)
\bezier{150}(30,5)(36,7)(38,10) \bezier{20}(30,15)(36,13)(38,10)
\bezier{40}(10,15)(20,17)(30,15)

\bezier{120}(2,10)(3,15)(7,20) \bezier{120}(13,35)(13,27)(7,20)
\put(13,35){\line(0,1){34.00}}

\bezier{120}(38,10)(37,15)(33,20) \bezier{120}(27,35)(27,27)(33,20)
\put(27,35){\line(0,1){34.00}}

\bezier{80}(15,33)(20,31)(25,33)

\bezier{80}(15,33)(13.5,33.7)(13,35)
\bezier{80}(25,33)(26.5,33.7)(27,35)
\bezier{7}(15,37)(13.5,36.3)(13,35)
\bezier{7}(25,37)(26.5,36.3)(27,35)

\bezier{20}(15,37)(20,39)(25,37)

\put(13.2,35){\line(1,0){13.6}} \put(17.6,32.6){\line(1,1){4.8}}

\put(2,10){\line(1,0){36.00}} \put(14.4,4.4){\line(1,1){11.20}}
\put(18.0,35.2){$\scriptstyle 1$} \put(21.5,67.5){$x$}
\put(26.5,21){$r$} \put(20,20){\vector(1,0){13.00}}
\put(37.0,50.0){(\emph{a})}
\put(85,70){\vector(0,-1){60.00}}
\bezier{120}(76,7.5)(85,4.5)(94,7.5)
\bezier{90}(76,7.5)(73.6,8.4)(73,10)
\bezier{90}(94,7.5)(96.4,8.4)(97,10)
\bezier{10}(76,12.5)(73.6,11.6)(73,10)
\bezier{10}(94,12.5)(96.4,11.6)(97,10)
\bezier{40}(76,12.5)(85,15.5)(94,12.5)
\put(73.0,10){\line(1,0){24.0}} \put(81.2,13.8){\line(1,-1){7.5}}

\bezier{120}(76,47.5)(85,44.5)(94,47.5)
\bezier{90}(76,47.5)(73.6,48.4)(73,50)
\bezier{90}(94,47.5)(96.4,48.4)(97,50)
\bezier{10}(76,52.5)(73.6,51.6)(73,50)
\bezier{10}(94,52.5)(96.4,51.6)(97,50)
\bezier{40}(76,52.5)(85,55.5)(94,52.5)

\put(73.0,50){\line(1,0){24.0}} \put(81.2,53.8){\line(1,-1){7.5}}
\put(73.0,50){\line(0,-1){40.0}} \put(97.0,50){\line(0,-1){40.0}}

\bezier{140}(73.0,50)(73,60)(78,70.5)
\bezier{140}(97.0,50)(97,60)(92,70.5)
\bezier{80}(80,68)(85,67)(90,68) \bezier{80}(80,68)(76,70)(80,72)
\bezier{80}(90,68)(94,70)(90,72) \bezier{80}(80,72)(85,73)(90,72)
\put(86.0,50.5){$\scriptstyle 1$} \put(86.5,30.5){$x$}
\put(90.0,61){$r$} \put(85,60){\vector(1,0){11.00}}
\put(105.0,50.0){(\emph{b})}
\end{picture}

\caption{\footnotesize  (a) $\D$ has eigenvalues, (b) $\D$ has not
eigenvalues } \label{Fig1}
\end{figure}

We consider the  distribution of resonances of $- \Delta_M$. For an analytic function $f$ on ${\C}$,
 let $\cN_r(f)$ be the number of zeros of $f$ counted with multiplicities in the disc of radius $r$ centered at the origin. By virtue of Theorem
\ref{T1}, i) and  the well-known result of Zwoski
\cite{Z87} for the distribution of resonances, we obtain

\begin{corollary}\lb{T3}
Let $q\in \cP_1$. Then
\[
\lb{Le} \cN_r(\p)={2r\/\pi}(1+o(1)) \ \ \qqq \as\qq r\to \iy.
\]

\end{corollary}

 Let $n_+(f)$  be the number of zeros of $f$ in
$\C_+$ counted with multiplicities.

\begin{corollary}
\lb{T2}
 Let $q\in \cP_1$. Then  the Jost function $\p(k )$ is represented as
\[
\lb{Hf} \p(k )=k^{n_0} \p^{(n_0)}(0)e^{ik }\lim_{r\to
\iy}\prod_{|k_n|\leq r, k_n\neq 0} \lt (1-{k \/ k_n} \rt), \ \ \ k
\in \C ,
\]
uniformly on any compact subset of $\C$, where $n_o=n_o(\p) = 0 \ or
\ 1$. Moreover, the S-matrix  has the form $S(k)=e^{-2i\f_{sc}(k)}$,
where $\f_{sc} $ is given by
\[
\lb{1.12}
\begin{aligned}
 \f_{sc}(k)=-\pi \lt(n_++{n_o\/2}\rt)+ \int_0^k
\f_{sc}'(t)dt,\ \ \
\\
\ \ \ \ \f_{sc}'(k )=1+\sum_{n=1,k_n\neq 0}^\iy {\Im k_n\/ |k
-k_n|^2} ,\ \ \ k>0,
\end{aligned}
\]
uniformly on any compact subset of $\R_+$, where $n_+=n_+(\p)$. Moreover, the following trace formula holds true:
\[
\lb{tr} -2k\Tr \rt(R(k)-R_0(k)\rt)={n_o\/k}+i+\lim_{r\to \iy}
\sum_{k_n\ne 0, |k_n|<r} {1\/k -k_n},
\]

\end{corollary}

\subsection{Inverse problem }
We sometimes write $\p(k,q), k_{n} (q), \cdots$ instead of $\p(k) ,k_{n}, \cdots$, when several potentials are dealt with. We introduce the
Sobolev space of real functions
\begin{equation}
\begin{aligned}
\lb{dWH} &W_1^0 =\Big\{q \in L_{real}^2(0,1)\ ;\ q'\in L^2(0,1), \
q(0)=q(1)=0\Big\},
\end{aligned}
\end{equation}
equipped with norm
$\|q\|^2_{W_1^0}=\|q'\|^2=\int_0^1|q'(x)|^2dx$. The main result of this paper is the following inverse problem.

\begin{theorem}
\lb{T4}
 i) Let $\p_j$ be the Jost function  for  $ q_j\in
W_1^0\cap \cP_1, j=1,2$. If $\p_1=\p_2$, then $q_1=q_2$.

\no ii)  Let $S_j$ be the  S-matrix for $
q_j\in W_1^0\cap \cP_1, j=1,2$. If $S_1=S_2$, then
$q_1=q_2$.

\no iii) Any $q\in W_1^0\cap \cP_1$ is uniquely determined by
its eigenvalues and resonances. Moreover, there exists an algorithm
to recover $q$ from its eigenvalues and resonances.
\end{theorem}

\subsection{Historical review}
There is an abundance of works devoted to the spectral theory and
inverse problems for the surface of revolution from the view points
of classical inverse Strum-Liouville theory, integrable systems,
micro-local analysis, see \cite{AA07}, \cite{E98}
and references therein. For integrable systems associated with
surfaces of revolution, see e.g. \cite{KT96}, \cite{T97},
 \cite{SW03} and references therein.

Isozaki-Korotyaev \cite{IK17}, \cite{IK17x}  solved the inverse
spectral problem for rotationally symmetric manifolds (finite
perturbed cylinders), which includes a class of surfaces of
revolution, by giving an analytic isomorphism from the space of
spectral data onto the space of functions describing the radius of
rotation. An analogue of the Minkowski problem was also solved. In
another paper \cite{IK19} Isozaki-Korotyaev  studied inverse problems
for Laplacian  on the torus. Moreover, they obtained stability
estimates: the spectral data in terms of the profile (the radius of
the rotation) and conversely, the profile in term of the spectral
data.

 In our paper we use inverse resonance theory for  Schr\"odinger operators  with
compactly supported potentials on the half-line \cite{K04}. We use
also results on perturbed Riccati mappings from \cite{K02},
\cite{K03}, \cite{IK17x}.

A lot of papers are devoted to resonances of  one-dimensional
Schr\"odinger operators. See Froese \cite{F97}, Hitrik \cite{H99},
Korotyaev \cite{K04}, Simon \cite{S00}, Zworski \cite{Z87} and
references therein.
 Inverse problems (uniqueness, reconstruction, characterization) in terms of resonances were solved by Korotyaev for a
Schr\"odinger operator with a compactly supported potential on the
real line \cite{K05} and the half-line \cite{K04}. See also Zworski
\cite{Z02}, Brown-Knowles-Weikard \cite{BKW03} concerning the
uniqueness. The resonances for  one-dimensional operators
$-{d^2\/dx^2}+V_\pi+V$, where $V_\pi$ is periodic and $V$ is a
compactly supported potential were considered  by Firsova
\cite{F84}, Korotyaev \cite{K11}, Korotyaev-Schmidt \cite{KS12}.
Christiansen \cite{C06} considered resonances for steplike
potentials. The "local resonance" stability problems were considered
in \cite{Ko04}, \cite{MSW10}.

 Resonances for specific cases of surfaces of revolution are discussed in
\cite{C04}, \cite{DKK15}, \cite{DH13}, \cite{D16}. As far as the authors know, the results in this paper about resonances for Laplacian on  surfaces of
revolution for the case \er{1}-\er{dP} are new.

\section {Preliminaries}
\setcounter{equation}{0}

\medskip

\subsection {Entire functions}
We recall some well-known facts about entire functions (see
\cite{Ko88}). An entire function $f(k)$ is said to be of
$exponential$ $ type$ if there is a constant $\a$ such that
$|f(k)|\leq $ const. $e^{\a |k|}$ everywhere on $\C$. An enitire function $f$ is
said to belong to the Cartwright class $E_{Cart},$ if $f(k)$ is of exponential type, and the following conditions hold true:
\[
\lb{14x}
\begin{aligned}
& \int _{\R}{\log ^+|f(x)|dx\/ 1+x^2}<\iy,
\\
& \r_+(f)=0,\ \ \ \r_-(f)=2,\ \ \ {\rm where}\ \ \ \r_{\pm}(f) =
\lim \sup_{t\to \iy} {\log |f(\pm it)|\/t}.
\end{aligned}
\]
Denote by $(z_n)_{n=1}^{\iy} $ the sequence of  zeros $\neq 0$ of $f$
(counted with multiplicities), arranged so that $0<|k_1|\leq |k_2|\leq
\cdots$. Then $f\in E_{Cart}$ has the Hadamard factorization
\[
\lb{2.7}
 f(k)=z^mCe^{iz}\lim_{R\to +\iy}\prod_{|k_n|\leq
R}\lt(1-{k\/ k_n}\rt),\ \ \ \ C={f^{(m)}(0)\/m!},
\]
for some integer $m\ge 0$, where the product converges uniformly on any compact set. Moreover, we have
\[
\sum_1^\iy {|\Im k_n|\/|k_n|^2}<\iy . \lb{2.8}
\]
Hence we obtain
\[
{f'(k)\/ f(k)}=i+{m\/ k}+\lim_{R\to \iy}\sum_{|k_n|\leq R}{1\/k-k_n}
\lb{2.9}
\]
uniformly on any compact subset of $\C\sm\big(\{0\}\cup\bigcup
\{k_n\}\big)$. In this paper, we make use of a simple class  of
functions in $E_{Cart}$. Let $g\in \cP_0$. Then
$1+\wh g$ belongs to $E_{Cart}$:
\[
\lb{c1x} 1+\wh g\in E_{Cart}, \qqq {\rm where} \qq \wh g(k)=\int_0^1
g(t)e^{2ikt}.
\]

Let us recall the Levinson Theorem (see \cite{Ko88}):
  {\it  Let $f\in E_{Cart}$. Then}
\[
\lb{c1} \cN (r,f)={2\/ \pi }r+o(r)\qqq \as \qq r\to \iy.
\]

\subsection {Unitary transformations.}
Recall that the Laplacian on $(M, g)$ acting in $L^2(M)=\os_{\n\ge
1} \mL_\n^2(M)$  is unitarily equivalent to a direct sum \er{sum} of
one-dimensional  operators $\D_\n$. Here the operator $\D_\n$ acts
in the space $L^2(\R_+,r^m(x)dx)$ under the  Dirichlet boundary
condition and is given by
\[
\begin{aligned}
\lb{3} & - \D_\n f=-{1\/r^m}(r^m f')' +u_\n f ,\qqq f(0)=0,
\\
& u_\n={E_\n\/r^2}=u_{\n,0}e^{-{4\/m}Q}, \qq u_{\n,0}=u_n(1),\qqq
 r=r_oe^{{2\/m}Q}, \qqq Q=-\int_x^1 q(t)dt.
\end{aligned}
\]
The function $\r=r^{m\/2}>0$ is called an impedance function.
 We define a unitary transformation $\mU$ by
 $$
 \mU:
L^2(R_+,\r dx)\to L^2(R_+,dx),\qqq \mU f= \r f,\qqq \r=r^{m\/2}=\r_0
e^{Q}.
$$
 We transform the operator $- \D_\n$ into a Schr\"odinger oprtator $H_\n$ by
\[
\lb{5}
\begin{aligned}
\mU (-\D_\n) \mU^{-1}=
 -\r^{-1}\pa_x \r^2\pa_x \r^{-1}+u_\n=\cD^*\cD+u=H_\n+u_{\n,0},
\end{aligned}
\]
 where $H_\n$ acts on $ L^2(R_+)$   and is given by
\[
\begin{aligned}
\lb{6}  H_\n =-{d^2\/dx^2}+p_\n, \qqq
p_\n=q'+q^2+u_\n-u_{_\n,0},\qqq
\end{aligned}
\]
and where the potential $p_\n\in L^1(\R_+)$ has a compact support $\ss
[0,1]$. Here we have used the identities
\[
\begin{aligned}
\lb{a7} & \cD=\r \pa_x \ \r^{-1}= \pa_x -q,\qqq \cD^*= \big(\r\
\pa_x \r^{-1}\big)^*=-\pa_x-q,
\\
& \cD^*\cD=-(\pa_x +q)(\pa_x -q)=-\pa_x^2+q'+q^2.
\end{aligned}
\]


The following support property of the non-linear  mappings $q \to
p_{\nu}$ in \er{6} plays a key role.

\begin{lemma}
\lb{Tqu}

 Let $q,q'\in L^1(\O)$ for
$\O=(1-\t,1)$ for some $\t
>0$ and $q(1)=0$. Assume that $q$ satisfies on $\O$ the following equation
\[
\lb{qc} q'+q^2+u-u_0=0,
\]
where $u=u_0 e^{\b Q}$ with $\b, u_0>0$. Then $q=0$ on
$\o=(1-\ve^2,1)$ for $\ve
>0$ small enough.

\end{lemma}
{\bf Proof.}
Define the norm in $L^r(\o), r\ge 1$ by  $\|q\|_r^r=\int_\o
|q(t)|^rdt$. We have for $x\in \o$
\[
\lb{q1}
   q(x)=-\int_x^1 q'(t)dt,\qqq |q(x)| \le \|q'\|_r \ve,\qqq
\|q\|\le \ve\|q'\|_1 ,
\]
where $\|q\|=\|q\|_2$. Let $q\ne 0$. Then for $\ve >0$ small enough
this yields
\[
\begin{aligned}
& \|q\|\le 1,\qqq    |Q(x)|\le \int_x^1 |q(t)|dt\le \ve \|q\|\le
{\ve} ,
\\
& |e^{\b Q(x)}-1| ={1\/\b}\rt|\int_x^1 q(t)e^{\b Q(t)}dt\rt|\le
{e^{\b \ve}\/\b}\int_x^1 |q(t)|dt \le {\ve\/\b}\|q\| e^{\b \ve},
\end{aligned}
\]
for  all $x\in \o$. Using these estimates for \er{qc},  we
have
\[
\begin{aligned}
& q(x)=-\int_x^1(q^2+u-u_0)dt,
\\
& |q(x)|\le \int_x^1 (|q(t)|^2+ u_0|e^{\b Q(x)}-1|)dt\le
\|q\|^2+{\ve \|q\|_\o u_0\/\b}e^{\b \ve}\le\|q\| C,
\end{aligned}
\]
 where  $C=1+{\ve  u_0\/\b}e^{\b \ve}$.
This  gives  $ \|q\|\le \ve\|q\| C$, which yields a contradiction,
since $\ve C< 1$   for $\ve >0$ small enough. \BBox


\section {Proof of main theorems}
\setcounter{equation}{0}

 \subsection {Eigenvalues and resonances}


 {\bf Proof of Theorem \ref{T1}.}
i) Let $q\in \cP_1$. We show that $p=p_\n$ given by \er{dP1} belongs
to $\cP_0$.   Assume that $p\notin \cP_0$ and $p=0$ on $\O=(1-\t,1)$
for some $\t >0$. Then by Lemma \ref{Tqu}, we have $q=0$ on
$\O=(1-\t_1,1)$ for some small $\t_1 >0$, which gives the
contradiction, since $q\in \cP_1$. Thus we have $p\in \cP_0$.

ii) The Schr\"odiger operator $H_\n=H_0+U_\n$, where $
H_0=-{d^2\/dx^2}+X_0\ge 0, \qq X_0=q'+q^2,
$
and the potential $U_\n$ has the form:
\[
\lb{U1}
\begin{aligned}
 & U_\n=u_\n-u_{_\n,0}=E_\n\rt({1\/r^2}- {1\/r_o^2} \rt)
=E_\n{r_o^2-r^2\/r_o^2r^2},\qq \supp U_\n\ss [0,1].
\end{aligned}
\]
 Recall that $u=u_{_\n,0}e^{-{4 Q\/m}}$ and $ Q=-\int_x^1 q(t)dt$.
Roughly speaking we have two cases:

1) $r<r_o$ on the interval $[0,1]$ and $U_\n\ge 0$, recall that
$r_o=r(1)$,

2) $r>r_o$ on the interval $[0,1]$ and $U_\n<0$.

Consider the case $r(x)\le r_{o}$ for all $x\in [0,1]$. Then the
potential $U_\n\ge 0$ and the operator $H_\n\ge H_0\ge 0$ and $H_\n$
has no eigenvalues. Note that $H_0\ge 0$ and the operator $H_0$ has
only absolutely continuous spectrum.


iii) Consider the second case $r(x)> r_{o}=1$ for all $x\in (0,1)$.
We take $r(x)=1+{(x-1)^2\/2}$ for $x\in [0,1]$. Then we obtain $1\le
r\le {3\/2}$ and $r'=x-1,\ r''=1, $ and the identity $\r=r^{m\/2}$
gives
\[
\lb{bx}
\begin{aligned}
& q={\r'\/\r}={m\/2}{r'\/r}={m\/2}{(x-1)\/r},
\\
&
q'={m\/2}{r''\/r}-{m\/2}{{r'}^2\/r^2}={m\/2r}-{2\/m}q^2={m\/2}-{r\/m}q^2-{2\/m}q^2=
{m\/2}-{r+2\/m}q^2,
\\
&  X_0=q'+q^2={m\/2}-{r+2-m\/m}q^2
\end{aligned}
\]
and
$$
U_\n=E_\n {1-r^2\/r^2}=-E_\n {(1+r)(x-1)^2\/2r^2}=-2E_\n
{(1+r)\/m^2}q^2.
$$
Let $\c_{\o}$ be the characteristic function of the interval $\omega
= (0,{1\/2})$. We have $q^2\ge {1\/4}\c_{\o}$. This and $p=X_0+U_\n$
yield
\[
\lb{bxx}
\begin{aligned}
p={m\/2}-{r+2-m\/m}q^2-2E_\n {(1+r)\/m^2}q^2={m\/2}-v_\n,
\\
v_\n=2E_\n q^2 \wt v_\n, \qqq \wt v_\n={1+r\/m^2}+{r+2-m\/m 2E_\n},
\qqq {1\/m^2}\le \wt v_\n \le {6\/m^2}.
\end{aligned}
\]
For $E_\n$ large enough we obtain $v_\n=2E_\n q^2 \wt v_\n\ge
v_\n^0:={E_\n\/2m^2}\c_{\o}$. Define the operators $h_\n, h_\n^0$ on
$\R_+$ by
$$
\textstyle  h_\n y=-y''+({m\/2}-v_\n) y, \qqq  h_\n^0
y=-y''+({m\/2}y-v_\n^0) y,\qq y(0)=0.
$$
It is well-known that the operator $h_\n^0$ has a large number of
negative eigenvalues for $E_\n$ large enough. Then the operator
$H_\n$ has also a large number of eigenvalues for $E_\n$ large
enough, since $H_{\n+1}\le H_\n\le h_\n\le h_\n^0$. This yields
iii).
 \BBox

{\bf Proof of Corollary \ref{T3}.} Theorem \ref{T1}, i) gives that
$p\in \cP_0$. In this case Zworski \cite{Z87} proved  the
asymptotics \er{Le} holds true. \BBox


{\bf Proof of Corollary  \ref{T2}.} Theorem \ref{T1}, i) gives that
$p\in \cP_0$. Thus all  results of Corollary   \ref{T2} follow from
Theorem \ref{Ta1} and \er{trx}. \BBox

 \subsection {Inverse problem}
We introduce the Sobolev spaces of real functions
\begin{equation}
\begin{aligned}
\lb{dH} \mH_0 =\rt\{q\in L_{real}^2(0,1): \  \int _0^1q(x)dx=0\rt\}
\qq \mH_\a =\rt\{q, q^{(\a)} \in \mH_0, \rt\}, \qqq \a \geq 0,
\end{aligned}
\end{equation}
equipped with norm
$\|q\|^2_{\mH_\a}=\|q^{(\a)}\|^2=\int_0^1|q^{(\a)}(x)|^2dx. $

We rewrite \er{dP1} in the form, omitting $\n$ for simplicity,
\[
\lb{dP2}
\begin{aligned}
p=v+v_0, \qqq v=q'+q^2+u-u_0-v_0=q'+q^2+u-c_0,
\\
 v_0=\int_0^1(q'+q^2+u-u_0)dx=\int_0^1(q^2+u-u_0)dx,
 \\
c_0=u_0+v_0=\int_0^1(q^2+u)dx.
\end{aligned}
\]
Thus we can consider the mapping $V:W_1^0\to \mH_0$ given by
\[
\begin{aligned}
\lb{dV}
& v=V(q)=q'+q^2+u-c_0,\\
& u(Q)=u_0e^{-\b Q},\qq Q(x)=-\int_x^1q(t)dt,\qq
 c_0=\int_0^1(q^2+u)dx,\qq \b={4\/m},
\end{aligned}
\]
where   $u_0={E_\n\/r_o^2}$.
The following result on the  mapping
$q\to v=V(q)$ is proven in  \cite{IK17}.

\begin{theorem}
\lb{T3x} The mapping $V:W_{1}^0\to \mH_0$, given by \er{dV} is a
real analytic isomorphism between the Hilbert spaces $W_{1}^0$ and
$\mH_{0}$ and satisfies:
\[
\lb{t3e} \|q'\|^2\le  \|v\|^2 \le
\|q'\|^2+2\|q\|^3\|q'\|+C_*\|q\|^2e^{2\b \|q\|},
\]
where the constant $C_*=u_0(\b+1)(2+\b u_0)$.

\end{theorem}

We consider the inverse problem: to determine a coefficient $q\in
W_0^1$, if the resonances of the operator $\D_\n$ are
given. Due to Theorem  \ref{T3x} it is equivalent to know the
resonance of the operator $H_\n, p_\n\in \mH_0$.

{\bf Proof of Theorem \ref{T4}.} i)  Let $\p_j(\cdot)$ be the Jost
function  for  some $ q_j\in W_1^0\cap \cP, j=1,2$. Assume
that $\p_1=\p_2$. Then Theorem \ref{Ta1} gives that  $p_1=p_2$. Thus Theorem \ref{T3x} implies that
$q_1=q_2$.

\no ii)  Let $S_j$ be the  S-matrix for  some  $
q_j\in W_1^0\cap \cP, j=1,2$. Assume that $S_1=S_2$. Let
$\p_j$ be the associated Jost function.
Let $\cZ(\p)$ be the sequence of zeros of $\p$ counted with
multiplicity. From the identity $S_1(k)={\p_1(-k)\/\p_1(k)}$ and
 \er{zj} we have $\cZ(\p_1)=\cZ(\p_2)$. Then
from \er{1.11} we obtain $\p_1=\p_2$. Then, i) gives $p_1=p_2$.
Thus Theorem \ref{T3x} implies  $q_1=q_2$.

\no iii) Let $\cZ(\p)$ be the set of all eigenvalues and
resonances for the  Jost function $\p(k,q)$ for some  $q\in W_1^0$
and the corresponding potential $p$. Then we can recover the Jost
function $\p$ by the Hadamard factorization \er{Hf}. Thus using the
Jost function $\p$ and the Marchenko equation (see
\er{Ma1}-\er{Ma4}), we can recover the potential $p$. If we know the
potential $p$, then we can recover the unique function $q\in
W_1^0\cap \cP$ by the perturbed Riccati equation using Theorem
\ref{T3x}. \BBox


\

\bigskip


\section {Inverse resonance scattering on the half-line}
\setcounter{equation}{0}

\bigskip

\subsection{Resonance scattering}
We consider the Schr\"odinger operator on $L^2(\R_+)$ given by
\[
Hy=-y''+p(x)y,\qqq y(0)=0,
\]
 where the potential $p\in \cP_0$.
 It is well-known (see \cite{F59}) that $H$ has purely
absolutely continuous spectrum $[0,\iy)$ plus a finite number
$n_+\geq 0$ of bound states $E_1<...<E_{n_+}<0$. Introduce the Jost
solution $f_+(x,k )$ and the regular solution $\vp (x,k)$ of the
equation
\[
-y''+py=k ^2y,\ \ \ x\geq 0, \ \ \ k \in \C\sm \{0\} ,
\]
with the conditions $f_+(x,k )=e^{ixk },\ x\geq 1$ and $\vp (0,k)=0,
\ \vp '(0,k)=1$. Outside the support of $p$ they are linear
combinations of $e^{\pm ik x}$. We define the Jost function $\p$ by
\[
\p(k)=f_+(0,k )= 1+\int _0^1{\sin k x\/k}p(x)f_+(x,k )dx, \ \ \ k
\in \C .
\]
This $\p(k)$ has the well known  properties, which we summarize
here.

\begin{lemma}
\lb{TaJ} If $p\in \cP_0$, the Jost function $\p(k)$ for the
operator $-y''+py$ on the half-line  with condition $y(0) = 0$ has the form
\[
\lb{fa1} \p(k)=1+{\hat F(k)-\hat F(0)\/2ik},\qq \forall \ k\in \C,
\]
for some $F\in \cP_0$. Moreover, the function $\p$ satisfies

 1) $\p(k)=1+{O(1)\/|k|}$ as $|k|\to \iy, k\in \ol\C_+$, uniformly in
 $\arg k\in [0,\pi]$.

2) $\p(k)$ is real on $i\R_+$,

3)  $\p(k)$ has $n_+$ simple  zeros in $\C_+$, which belong to
$i\R_+$,

4) the function $\f_{sc}$ is odd $ \f_{sc}(-k )=-\f_{sc}(k ), \ k>0,
$ and is continuous on $(0,\iy )$

5) $\f_{sc}(k )=o(1)$ as $k\to \iy $ and
\[
\f_{sc}(\pm 0)=\mp \pi \lt (n_++{n_o\/2}\lt),\ \ \ \as\qq  n_o=n_o
(\p)=0, 1.
\]

6) Let $\cZ(f)$ be the set of zeros of an entire function $f$. The
zeros of $\p$ satisfies
\[
\lb{zj}
\cZ(\p)\cap \cZ(\p(-\cdot))=\ca \{0\} \qq & if \qq \p(0)=0, \\
                            \emptyset & if \qq  \p(0)\ne 0. \ac
\]

\end{lemma}

Introduce the set $\cJ$ of all possible Jost functions from
\cite{K04}.

\medskip
\no {\bf  Definition $\cJ$}. {\it By   $\cJ$ we mean the class of
all entire functions $f$  having the form
\[
\begin{aligned}
f(k )=1+{1\/ 2ik}\lt(\hat F(k)-\hat F(0)\rt), \ \ \ \ \  \ \ k\in
\C,
\end{aligned}
\]
where $\hat F(k)=\int_0^1 F(x)e^{2ixk }dx$ and $F\in \cP_0$. In
addition the sequence $(k_n)_{n=1}^\iy$ of zeros of $f$
satisfies:

\no i) only one $f(k)\neq 0$ for any $k \in\R\sm \{0\}$ and
$n_o(f)\leq 1$,

\no ii)   all zeros $k_1,..,k_{n_+}, n_+=n_+(f)$ of the function $f$
in $\C_+$ are simple, belong to $i\R_+$ and satisfy}
\[
\begin{aligned}
&  |k_1|>|k_2|>...>|k_{n_+}|>0,
  \\
&   
 (-1)^jf(-k_j) >0,\ \qqq j=1, ... , {n_+}.
\end{aligned}
\]

We recall the well-known Marchenko result \cite{M86} about the
inverse problems for our case $p\in \cP_0$. Define the functions
\[
\lb{Ma1}
  G_0(x)= {1\/ 2\pi }\int_{-\iy}^\iy (1-S(k ))e^{ixk}dk ,
\ \ \ \ \ G(x)=G_0(x)+\sum_{k_j\in\C_+} {e^{-x|k_j|}\/c_j} ,
\]
and the norming constants
\[
\lb{Ma2}
  c_j=\int_{\R_+}|\p (x,k_j)|^2dx<\iy , \ \ \ k_j\in \C_+, \qq j=1,.., n_+,
\]
where $\p (\cdot, k_j)\in L^2(\R_+), k_j\in \C_+$  is an
eigenfunction of $H$. For each $x>0$ the Marchenko equation
\[
\lb{Ma3} K(x,t)=-G(x+t)-\int_{x}^\iy G(t+s)K(x,s)ds, \ \ \ t\geq x,
\]
has the unique solution $K(x,t)$, which satisfies \er{K}, \er{K1}.
It is important that the potential $p$ has the form
\[
\lb{Ma4} p(x)=-2{d\/ dx}K(x,x),  \ \ \ x>0.
\]
Below we use the fact: if $G(x)=0$ for $x>2d$, then $p(x)=0$ for $x>d$.

 Now we consider the mapping
$p\to J(p)=\p(\cdot,p)\in \cJ$, i.e., for each $p\in \cP_0$ we have
the Jost function $J(p)=\p(k,p)$. We are now in a position to formulating the main result on the  inverse resonance scattering from \cite{K04} including
characterization.

\begin{theorem}
\lb{Ta1}
  i) The mapping $J: \cP_0\to \cJ$ given by
$p\to J(p)=\p(\cdot,p)$    is one-to-one and onto.

\no ii) For any $p\in \cP_0$ the Jost function $\p(k)=\p(k,p)$ is
given by
\[
\lb{1.11}
 \p(k)=k^{n_o} \p^{(n_o)}(0) e^{ik }\lim_{r\to
\iy}\prod_{|k_n|\leq r, k_n\neq 0} \lt (1-{k \/ k_n} \rt), \ \ \ k
\in \C ,   \ \ \ n_o={n_o}(\p) = 0, 1,
\]
where the product is uniformly convergent on any compact subset of $\C$.
Moreover, the potential $p$ is uniquely determined by the Marchenko
equations \er{Ma1}-\er{Ma4} (in terms of its bound states and
resonances), where
$S(k)=e^{-2i\f_{sc}(k)}$ and $\f_{sc} $ is represented as
\[
\lb{12}
\begin{aligned}
\f_{sc}(k)=-\pi \lt(n_+(\p)+{n_o\/2}\rt)+ \int_0^k \f_{sc}'(t)dt,\ \
\\
\ \ \ \ \f_{sc}'(k )=1+\sum_{n=1,k_n\neq 0}^\iy {\Im k_n\/ |k
-k_n|^2} ,\ \ \ k\geq 0,
\end{aligned}
\]
and the norming constant $c_j$  (see \er{Ma2}) is given by
\[
 \lb{1.13}
 c_j=(-1)^\s {e^{-2|k_j|}\/ 2|k_j|}\prod_{n\neq j, n\geq 1,k_n\neq 0}^\iy
 \lt ({1-{k_j \/ k_n}\/ 1+{k_j \/ k_n}} \rt),\ \ \ j=1, ..., n_+(\p).
\]
\end{theorem}

{\bf Remark.} 1)  The series \er{12} converges absolutely (see
(2.8)).

2) Moreover, the following trace formula holds true:
\[
\lb{trx} -2k\Tr \rt(R(k)-R_0(k)\rt)={{n_o}\/k}+i+\lim_{r\to \iy}
\sum_{k_n\ne 0, |k_n|<r} {1\/k -k_n},
\]
uniformly on compact subsets of $\C\sm\lt(\{0\}\cup\bigcup
\{k_n\}\rt)$.


 We divide the inverse problem for
the mapping  $k_{\bu} : p\to k_{\bu}(p)=(k_n(p))_1^\iy$ into three
parts:

\medskip
\noindent
\no I. {\it Uniqueness}. Do the eigenvalues and the resonances
determine uniquely $p$?

\smallskip
\noindent
\no II. {\it Reconstruction}. Give an algorithm for recovering $p$
from the sequence  $k_{\bu}(p)=(k_n(p))_1^\iy$ only.

\smallskip
\noindent
\no III. {\it Characterization}. Give necessary and sufficient
conditions for the sequence  $k_{\bu}=(k_n)_1^\iy$ to be the
eigenvalues and the resonances for some $p\in \cP_0$.

\medskip
Theorem \ref{Ta1} gives the solution of I and II. Moreover, it gives
a characterization of $p\in \cP_0$ in terms of functions from $\cJ$
through the mapping $J : \cP_0 \to \cJ$. The proof of Theorem
\ref{Ta1} is based on the entire function theory and the Marchenko
results about the inverse scattering on the half-line.

\medskip

\setlength{\itemsep}{-\parskip} \footnotesize \no  {\bf
Acknowledgments.}   E. Korotyaev was supported by the RSF grant  No.
18-11-00032.  H.Isozaki is partially supported by Grants-in-Aid for
Scientific Research (S) 15H05740, and Grants-in-Aid for Scientific
Research (B) 16H03944, Japan Society for the Promotion of Science.

\end{document}